\documentclass[11pt]{amsart}
\usepackage[greek,english]{babel}
\usepackage[utf8x]{inputenc}
\usepackage[T1]{fontenc}
\usepackage{a4wide}
\usepackage{amssymb}
\usepackage{mathpazo}
\usepackage{graphicx}
\usepackage[dvipsnames]{xcolor}
\usepackage[colorlinks=true,allcolors=MidnightBlue,linktoc=page]{hyperref}
\usepackage{tikz-cd}
\usepackage{pict2e}

\everymath{\displaystyle}

\DeclareMathOperator{\SL}{SL}

\newcommand{\action}{\curvearrowright}

\DeclareMathOperator{\co}{co}
\DeclareMathOperator{\Sym}{Sym}

\DeclareMathOperator{\Ends}{Ends}
\DeclareMathOperator{\Stab}{Stab}

\DeclareMathOperator{\Aut}{Aut}

\DeclareMathOperator{\Br}{Br}

\DeclareMathOperator{\Id}{Id}

\DeclareMathOperator{\Homeo}{Homeo}
\DeclareMathOperator{\Diff}{Diff}

\newcommand{\R}{\mathbf{R}}                          
\newcommand{\C}{\mathbf{C}}        
\newcommand{\D}{\mathbf{D}}                     
\newcommand{\N}{\mathbf{N}}                          
\newcommand{\Z}{\mathbf{Z}}                          
\newcommand{\bS}{\mathbf{S}}   

\newtheorem{theorem}{Theorem}[section]
\newtheorem{lemma}[theorem]{Lemma}

\newtheorem{corollary}[theorem]{Corollary}

\newtheorem{proposition}[theorem]{Proposition}

\theoremstyle{definition}
\newtheorem{definition}[theorem]{Definition}

\newtheorem{remark}[theorem]{Remark}

\numberwithin{equation}{section}

\author{Bruno Duchesne}
\date{April 2021}
\thanks{Pierre Py emphasized the question of groups with Property (T) acting on the circle and made comments on a preliminary version. Michael Megrelishvili helped me to spot WAP triviality of $\Homeo^+(\bS^1)$. Peter Haïssinsky helped me to find my way in the literature about complex quadratic dynamics. Nicolas Monod and Andrés Navas made relevant comments on a preliminary version. Andrés Navas asked if the action can be smoothened. Víctor José García Garrido provided Figure~\ref{Garcia}.  Gianluca Basso found a mistake in the proof of minimality in Theorem~\ref{umf}. The mistake has been corrected using his ideas. I would like to thank all of them.}
\begin{document}
\title{A group with Property (T) acting on the circle}
\maketitle
\begin{abstract}
We exhibit a topological group $G$ with property (T) acting non-elementarily and continuously on the circle. This group is an uncountable totally disconnected closed subgroup of $\Homeo^+(\bS^1)$. It has a large unitary dual since it separates points. It comes from homeomorphisms of dendrites and a kaleidoscopic construction. Alternatively, it can be seen as the group of elements preserving some specific geodesic lamination of the hyperbolic disk.

We also prove that this action is unique up to conjugation and that it can't be smoothened in any way. Finally, we determine the universal minimal flow of the group $G$.  \end{abstract}
\section{Introduction}
The Zimmer program aims to classify actions of higher rank Lie groups and their lattices on compact manifolds. Originally Zimmer considered actions by diffeomorphisms but one can consider more generally actions by homeomorphisms. In this setting, the simplest compact manifold is the circle. First results in this direction were obtained by Witte for $\SL_n(\Z)$ for $n\geq3$ \cite{MR1198459}. He showed that any homomorphism to the group of orientation-preserving homeomorphisms $\Homeo^+(\bS^1)$ has finite image. Latter Ghys \cite{MR1703323,MR1876932} proved that lattices in  higher rank simple connected Lie groups have finite orbits.  Since such lattices  have property (T), it is natural to ask if these rigidity results could be a consequence of Property (T). See Question 2 in \cite{MR3966841}. For discrete groups acting by $C^{1+\alpha}$-diffeomorphisms (where $\alpha>1/2$), Navas proved that the image of any Property (T) group in $\Diff(\bS^1)$ is finite \cite{MR1951442}. 

Very recently Deroin and Hurtado proved that higher rank lattices of semi-simple Lie groups with finite center are not left-orderable \cite{bertr2020non}. As a consequence they proved that lattices in  higher rank simple connected Lie groups have finite image. The question of the existence of an infinite countable group with Property (T) acting faithfully on the circle is still open.\\

The aim of this note is to broaden the frame to uncountable subgroups of $\Homeo^+(\bS^1)$.  This group,  $\Homeo^+(\bS^1)$, has a natural topology, the compact-open topology and this topology is non-locally compact but still pleasant since it is Polish. It means that it is separable and metrizable by a complete metric. Since Property (T) is a topological property, it is natural to look for non-discrete such subgroups. We give an example of Property (T) closed subgroup $G$ of $\Homeo^+(\bS^1)$ which is non-elementary.

After Margulis version of the Tits alternative for groups acting on the circle, the right notion of elementarity for groups acting on the circle is the following one \cite[Théorème 3]{MR1797749}. 
\begin{definition} Let $\Gamma$ be a subgroup of $\Homeo(\bS^1)$. The subgroup $\Gamma$ is said to be \emph{elementary} if it preserves a probability measure.\end{definition} For example, all amenable subgroups of $\Homeo(\bS^1)$ are elementary.\\

Let us introduce the construction of the group $G$. Inspired by Burger-Mozes universal groups $U(F)$ for $d$-regular trees and a permutation group $F\leq \Sym([d])$ \cite{burger2000groups}, \emph{kaleidoscopic groups} were introduced in \cite{MR4017014}. One construct a group $\mathcal{K}(\Gamma)$ from a permutation group $\Gamma\leq\Sym([n])$ where $n\in\{3,4,\dots,\infty\}$ and $[n]=\{1,2,\dots,n\}$. The elements of this group $\mathcal{K}(\Gamma)$ act by homeomorphisms on the Wa\.zewski dendrite $D_n$ whose branch points have order $n$. Informally, dendrites can be thought as compact tress with a dense set of branch points. Dendrites and this kaleidoscopic construction are described in \S \ref{remind}.  Let $G$ be the Kaleidoscopic group $\mathcal{K}(A_3)$ associated with the alternating group $A_3$. This is a non-Archimedean Polish group, which means that it is isomorphic as topological group to a closed subgroup of the symmetric group $S_\infty$ with its pointwise convergence topology. Thus it is uncountable and totally disconnected. 

By a \emph{topological group embedding}, we mean a continuous injective group homomorphism $G\to H$ between two topological groups such that it is a homeomorphism on its image. In this case, we say that $G$ embeds topologically in $H$.

\begin{theorem}\label{action}The group $G$ embeds topologically in $\Homeo^+(\bS^1)$ as a non-elementary closed subgroup.
\end{theorem}

\begin{figure}[h]
\begin{center}
\includegraphics[width=.6\textwidth]{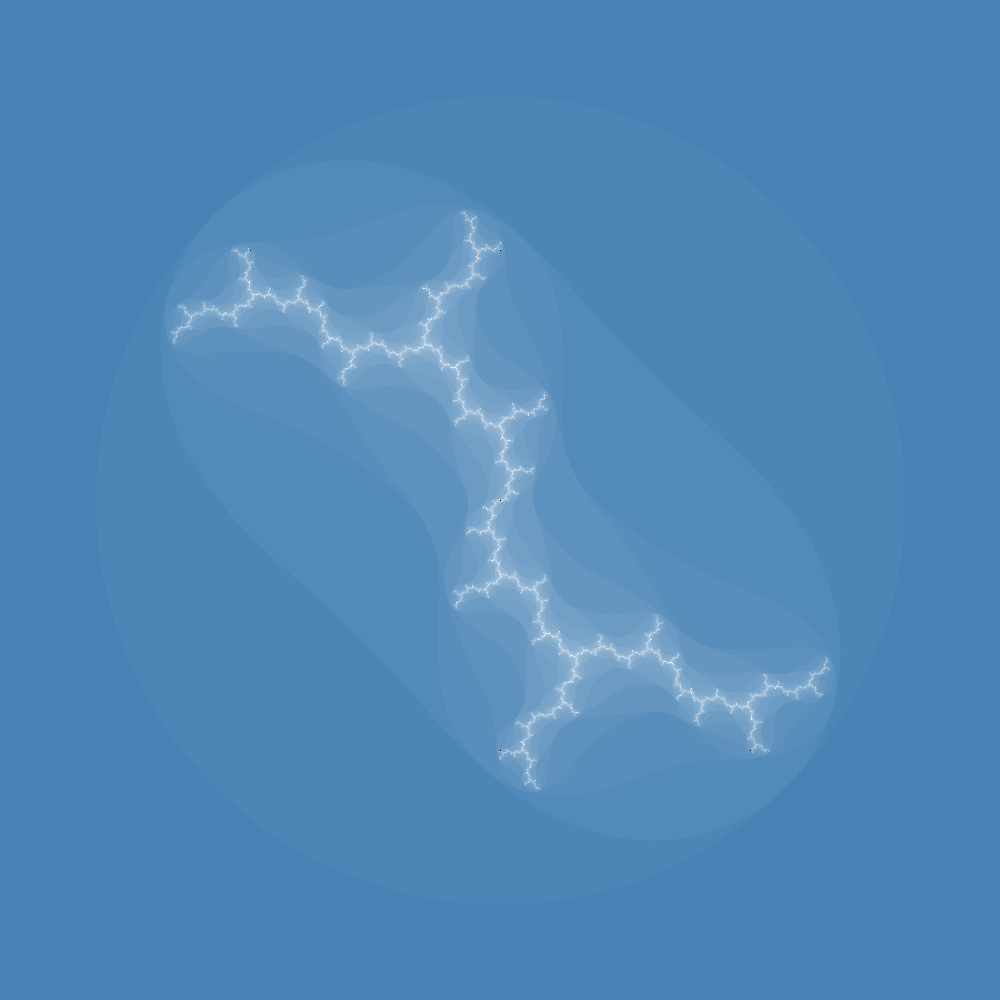}
\caption{The Wa\.zewski dendrite $D_3$ realized as the Julia set  of the complex polynomial $z^2+i$.}
\end{center}
\end{figure}

The main point of this theorem is the construction of an action $G\action \bS^1$. This action comes from a Carathéodory loop $\varphi\colon \bS^1\to D_3$ where $D_3$ is the  Wa\.zewski dendrite of order 3, realized as the Julia set of the polynomial $z^2+i$. The group $G=\mathcal{K}(A_3)$ acts on $D_3$ by definition and the action on $\bS^1$ is induced by $\varphi$. This action can be bonded to \cite[Problem 4.4]{MR1876932} where Ghys asks for a description of closed subgroups of $\Homeo^+(\bS^1)$ with only one orbit (see \cite{MR2255499} for an answer to this problem).

\begin{proposition}\label{orbits} The group $G$ has exactly 3 orbits in $\bS^1$.
\end{proposition}
A group action $\Gamma\action X$ is \emph{oligomorphic} if for any finite $n$, the diagonal action $G\action X^n$ has finitely many orbit. The action $G\action \bS^1$ is not transitive but it is very homogeneous in some sense.

\begin{proposition}\label{oligomorphy} The action $G\action\bS^1$ is oligomorphic.\end{proposition}

Let us recall that a topological group $G$ has \emph{Property (T)} if there is  $K\subset G$ compact and $\varepsilon>0$ such that for any continuous unitary representation $\pi\colon G\to \mathrm{U}(\mathcal{H})$, where $\mathcal{H}$ is a Hilbert space, the existence of a unit vector $v$ such that for all $g\in K$, $||\pi(g)(v)-v||<\varepsilon$ implies that $\pi$ has an invariant unit vector.

\begin{proposition}\label{T} The group $G$ has Property (T).\end{proposition}

One can actually see the group $G$ as the stabilizer of a lamination without direct reference to dendrites. Let us recall that  a \emph{laminational equivalence relation} is a closed equivalence relation on $\bS^1$ with finite equivalence classes and such that any two distinct classes have disjoint convex hulls in the disk. Let $\sim$ be such  a laminational equivalence relation. The automorphism group of $\sim$ is the subgroup of elements $g\in\Homeo^+(\bS^1)$ such that for all $z_1,z_2\in\bS^1$, $z_1\sim z_2\iff g(z_1)\sim g(z_2)$.  To a laminational equivalence relation, there is an associated \emph{geodesic lamination} of the hyperbolic disk. This is the collection of geodesics between points in the same equivalence class in the circle at infinity.


For the Carathéodory loop $\varphi$ used above, there is an associated laminational equivalence relation $\sim_\varphi$ defined by $z_1\sim z_2\iff \varphi(z_1)=\varphi(z_2)$.

\begin{proposition}\label{lamination}The group $G$ is the automorphism group of the lamination $\sim_\varphi$.\end{proposition}

\begin{remark} It is not difficult to find a Property (T) group that acts non-elementarily on $\bS^1$. Merely, because the group $\Homeo^+(\bS^1)$ itself has Property (T). This follows from a much stronger result due to Glasner and Megrelishvili. They proved that $\Homeo^+(\bS^1)$ is WAP-trivial, that is any weakly almost periodic function on $\Homeo^+(\bS^1)$ is constant \cite[Remarks 7.5.2]{MR3821724}. Since matrix coefficients are WAP functions, this implies that any unitary representations of $\Homeo^+(\bS^1)$ is trivial and thus this group has Property (T) for the trivial reason that it lacks unitary representations. This also explains the difficulty to answer negatively the question of the existence of countable Property (T) groups acting on the circle  since there is no unitary representations that exist for all such groups.
\end{remark}

\begin{remark}
For the group $G$, there is a wealth of unitary representations since it is non-Archimedean and thus unitary representations separate points \cite[Theorem 1.1]{MR2929072}. Actually, $G$ is a closed subgroup of $S_\infty=\Sym(\N)$ and the representation by permutation of the coordinates on $\ell^2(\N)$ is continuous and faithful.
\end{remark}

\begin{remark} One can actually construct infinitely many non-isomorphic groups with Property (T) embedded  in $\Homeo^+(\bS^1)$. Similar arguments show that the kaleidoscopic groups $\mathcal{K}(C_n)$ where $C_n$ is the cyclic group of order $n$ is oligomorphic and thus has Property (T), and acts faithfully continuously by orientation preserving homeomorphisms on $\bS^1$.\end{remark}

\begin{remark} This let open the question to determine whether a countable group with Property (T) may act non-elementarily (equivalently without a finite orbit in this case \cite[Hint of Exercise 5.2.12]{MR2809110}) on $\bS^1$. This question is reminiscent of the analogous question for groups with Property (T) acting on dendrites. We know such examples among non-Archimedean groups but no such countable groups are known \cite{DuchesneMonod,MR3894039}. 

Let us observe that the group $\SL_2(\Z)\ltimes \Z^2$ has relative property (T) with respect to $\Z^2$ and acts non-elementarily on the circle, see \cite[Example 2.106]{MR2327361} and \cite[Example 5.2.31]{MR2809110}. Moreover, this action can’t be smoothened to a $C^1$-action, see \cite[Remark 5.2.33]{MR2809110} and \cite{MR2602845}.\end{remark}

\begin{remark} For locally compact second countable groups, Property (T) is equivalent to Property FH, the fixed point property for continuous isometric actions on Hilbert spaces. This is no more the case in general \cite{MR2415834} (but Property (T) still implies Property FH). Rosendal proved that $\Homeo(\bS^1)$ has strong uncountable cofinality and thus any isometric action on a metric space has bounded orbits \cite{MR2503307}. Cornulier proved the same result for  $\Homeo(\bS^n)$, $n\geq 1$ \cite[Appendix]{MR2207794}. Thanks to the center lemma, Property FH follows from this property as well.\end{remark}


To conclude, one can ask if the group $G$ admits other continuous actions on the circle. In particular, can $G$ act on $\bS^1$ by $C^1$-diffeomorphisms or Hölder homeomorphisms?  We prove that the action is unique and thus, the action cannot be smoothened in any way. 

\begin{theorem}\label{unique} The action $G\action \bS^1$ described above is the unique minimal continuous action of $G$ on the circle up to conjugation by an element of $\Homeo^+(\bS^1)$.
\end{theorem}

As a consequence, we get that we cannot improve the regularity of the action. Let us recall that a \emph{modulus of continuity} is a map $\omega\colon \R^+\to\R^+$ such that $\omega(0)=0$ and $\omega$ is continuous at 0. For a modulus of continuity $\omega$,  a map between metric spaces $(X,d)$ and $(Y,d)$ has \emph{continuity type} $\omega$ if  there is $\lambda>0$ such that $\forall x,x’\in X$, $d(f(x),f(x’))<\lambda \omega(d(x,x’))$.  

For example, if $\omega(r)=r^\alpha$ with $\alpha>0$, maps with continuity type $\omega$ coincide with $\alpha$-Hölder maps.

\begin{theorem}\label{no smooth} For any modulus of continuity $\omega$ and any minimal continuous action $G\action\bS^1$, there is $g\in G$ such that the image of $g$ in $\Homeo^+(\bS^1)$ does not have continuity type $\omega$.\end{theorem}

While showing that the action $G\action \bS^1$ is unique up to homeomorphisms, we are very close to identify the \emph{universal minimal flow} of $G$. This is the unique (up to homeomorphism) minimal compact continuous $G$-space, i.e., minimal $G$-flow, such that any other minimal $G$-flow is a quotient of this one. 

Let $\xi$ be some end point of the dendrite $D_3$ and let $G_\xi$ be its stabilizer. We endow the quotient space $G/G_\xi$ with the uniform structure induced by the left uniform structure on $G$. A basis of entourages is given by 
$$U_V=\{(gG_\xi,hG_\xi),\ h\in vgG_\xi, \textrm{with}\ v\in V\}$$ where $V<G$ is the pointwise stabilizer of a finite set of branch points in the dendrite. We denote by $\widehat{G/G_\xi}$ the completion of this uniform structure.

\begin{theorem}\label{umf}The universal minimal flow of $G$ is $\widehat{G/G_\xi}$.
\end{theorem} 

This theorem follows essentially from the fact that $G_\xi$ is \emph{extremely amenable}, that is any continuous action on compact space has a fixed point and the understanding of the $G$-map $\widehat{G/G_\xi}\to D_3$. Since $G\action \widehat{G/G_\xi}$ is strongly proximal, we deduce the following.

\begin{corollary} Any minimal $G$-flow is strongly proximal.\end{corollary}

\noindent\textbf{Organization of the note.} The second section gathers the needed material about dendrites, kaleidoscopic groups and Carathéodory loops. The third section contains the construction of the action and the proof of the first statements presented in the introduction. The fourth section deals about essential uniqueness of this action and the last one, considers the side question of the determination of universal minimal flow of the group $G$.

\section{Wa\.zewski dendrites and Carathéodory loops}\label{remind}
\subsection{Wa\.zewski dendrites}
A \emph{dendrite} is a metrizable connected compact space that is locally connected and such any two points are connected by a unique arc. We refer to \cite{Nadler} for a general background. Dendrites can be thought as "compact topological trees." Let $X$ be a dendrite non-reduced to a point and let $x\in X$. A \emph{branch around} $x$ is a connected component of $X\setminus\{x\}$. The \emph{order} of $x$ is the number of branches around it. If this number is 1 then $x$ is an \emph{end point}, if it is 2 then $x$ is a \emph{regular point} and otherwise $x$ is a \emph{branch point}. We denote by $\Br(X)$ the set of branch points and for $x\neq y$, $B_x(y)$ is the branch around $x$ containing $y$. We also denote by $[x,y]$ the unique arc (injective continuous image of $[0,1]$) between $x$ and $y$ and by $(x,y)$, the open arc, that is $[x,y]\setminus\{x,y\}$. Let $D(x,y)$ be the connected component of $X\setminus\{x,y\}$ that contains $(x,y)$. Equivalently, this is $B_x(y)\cap B_y(x)$. For three points $x,y,z$ in a dendrite, the intersection $[x,y]\cap [y,z]\cap[z,x]$ is reduced to point that is called the \emph{center} of $x,y,z$.

 In the 1920s, Wa\.zewski introduced a universal dendrite that contains all the others \cite{Wazewski23}. This dendrite is characterized by the following properties: branch points are dense and their order is infinite (necessarily countable). Wa\.zewski construction can be generalized to finite cardinals and we denote by $D_n$ (for $n\in\{3,4,\dots,\infty\}$) the unique dendrite such that branch points are dense and have order $n$.

Since the original Wa\.zewski dendrite is constructed as a subset of the plane, all dendrites can be identified with a subset of the plane as well. Here we give two such identifications of the dendrite $D_3$. Let us start with a classical example coming from one-dimensional complex dynamics. For a complex polynomial $P$, the \emph{Julia set} of $P$ is the subset of the plane where the dynamic of $P$ is chaotic, more precisely, the subset of points $z\in\C$ such that the family of polynomials $P^n$ is not normal on neighborhoods of $z$. See \cite{MR2193309} for a classical reference about complex dynamics. The following result seems to be well known and it is claimed in \cite{bonk2018continuum} without proof nor references.

\begin{proposition} The Julia set of the polynomial $z^2+i$ is homeomorphic to $D_3$.\end{proposition}

\begin{proof} It is well known that this Julia set is a dendrite. See for example \cite[Example 1.1.5]{MR1747010} or \cite[Theorem V.4.2]{MR1230383}. Thanks to \cite{MR2508255}, any branch point is preperiodic or precritical and the only critical point is preperiodic. This point is $0$ and it is mapped successively to $i$, $-1+i$  and $-i$. Then we have a 2-cocycle: $-i\leftrightarrow -1+i$. The fixed point $\alpha_i=\frac{1+\sqrt{1-4i}}{2}$ (where we use the principal branch of the square function) of this polynomial is the landing point of the external rays $1/7,2/7$ and $4/7$ \cite[Figure 1]{MR1209913}. So this point has order 3 (see Lemma~\ref{branches} below). It follows from Thurston laminations technics that any branch point contains $\alpha_i$ in its positive orbit and thus has the same order as $\alpha_i$, that is 3. Since preimages of $1/7,2/7$ and $4/7$ by multiplication by 2 in $\R/\Z$ are dense, it follows that branch points in the dendrite are dense and this Julia set is homeomorphic to $D_3$. \end{proof}
\begin{figure}\label{Garcia}
\begin{center}
\includegraphics[width=.6\textwidth]{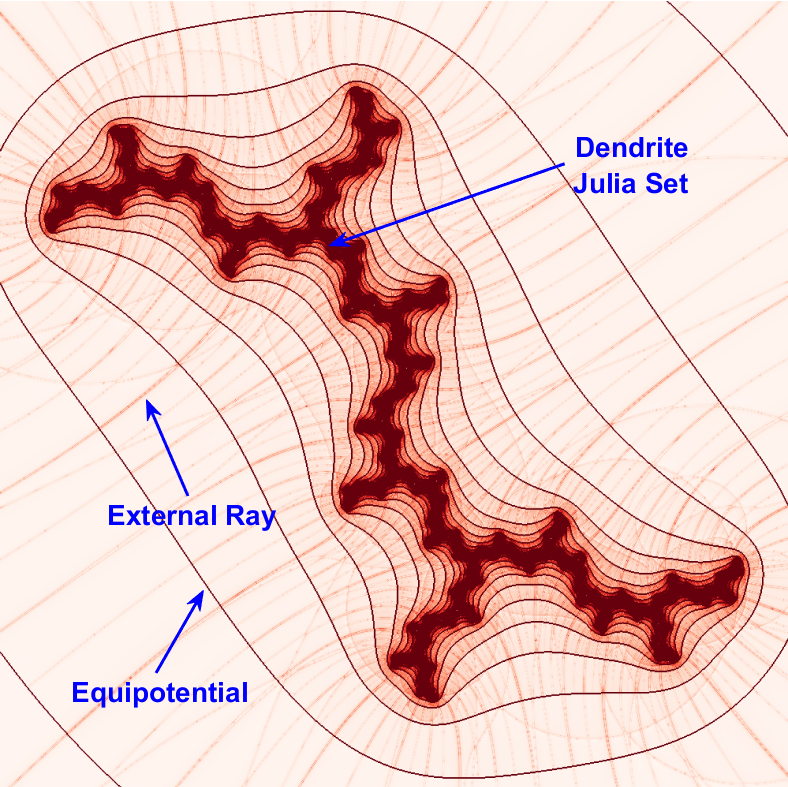}
\caption{The dendrite $D_3$ as Julia set of the polynomial $z^2+i$ with external rays and equipotential curves. Image kindly provided by Víctor José García Garrido. See \cite{GARCIAGARRIDO2020105417}.}
\end{center}
\end{figure}

We denote by $\widehat{\C}$ the Riemann sphere $\C\cup\{\infty\}$ with its conformal structure. The classical Riemann mapping theorem states that any simply connected domain $U$ of $\widehat{\C}$ with non-empty complement is the image of a conformal homeomorphism $\psi\colon \widehat{\C}\setminus\overline{\D}\to U$ where $\D$ is the open unit disk. Moreover if $\infty\in U$, the map $\psi$ is the unique such map under the conditions $\psi(\infty)=\infty$ and $\psi’(\infty)=1$. Let us identify $D_3$ with some subset of $\C$ and let $U$ be its complement. Since $U$ is simply connected and $D_3$ is locally connected, Carathéodory theorem (\cite[Theorem 17.4]{MR2193309})  shows that $\psi$ extends to a continuous map $\overline{\psi}\colon\widehat{\C}\setminus \D\to \overline{U}=\overline{\C}$. The restriction of this extension gives a continuous surjective map $\varphi\colon\bS^1\to D_3$ that is the called the \emph{Carathéodory loop} of $D_3$.  We refer to \cite[\S17]{MR2193309} for details. The images by $\overline{\psi}$ of rays $\{re^{i\theta},\ r\geq 1\}$ are called \emph{external rays}.

In the particular case of the Julia set of $z^2+i$, the Carathéodory loop semi conjugates the action of the polynomial maps $z\mapsto z^2$ and $z\mapsto z^2+i$, that is $\varphi(z^2)=\varphi(z)^2+i$ for all $z\in\bS^1$.

We fix a positive orientation on the circle. For $z_1\neq z_2$, we denote by $[z_1,z_2]_+$ the positive closed interval of points $z$ such that $(z_1,z,z_3)$ is positively oriented (or cyclically oriented). We denote by $(z_1,z_2)_+$ the open positive interval.

\begin{lemma}\label{branches} For $x\in D_3$

\begin{itemize}
\item $|\varphi^{-1}(\{x\})|=1\iff$ $x$ is an end point,
\item $|\varphi^{-1}(\{x\})|=2\iff$ $x$ is a regular point,
\item $|\varphi^{-1}(\{x\})|=3\iff$ $x$ is a branch point.
\end{itemize}

Moreover, if $x$ is a branch point of $D_3$ and $z_1,z_2,z_3$ are its preimages ordered cyclically on the circle then $\{\varphi((z_i,z_{i+1})_+)\}_{I=1,2,3}$ are the three branches around $x$. Similarly, if $x$ is a regular point and $z_1,z_2$ are its two preimages then $\varphi((z_1,z_2)_+)$ and $\varphi((z_2,z_1)_+)$ are the two branches around $x$.
\end{lemma}

\begin{proof} Let $z\neq z’\in\bS^1$ such that $\varphi(z)=\varphi(z’)=x$. By \cite[Theorem 17.4]{MR2193309} , there $z_1\in(z,z’)_+$ and $z_2\in(z,z’)_+$ such that $\varphi(z_i)\neq x$ and, moreover, these points lie in distinct connected components of $D_3\setminus\{x\}$  \cite[Lemma 17.5]{MR2193309}. 

For $x\in D_3$ and $y,z$ in two distinct branches around $x$ there are $z_1,z_2\in\bS^1$ such that $\varphi(z_1)=y$ and $\varphi(z_2)=z$. By connectedness of $(z_1,z_2)_+$ and continuity of $\varphi$, there is $z_3\in(z_1,z_2)_+$ such that $\varphi(z_3)=x$. Let $Z=\varphi^{-1}(\{x\})$. The last argument shows that preimages of branches around $x$ are in distinct connected components of $\bS^1\setminus Z$ and the former one shows that any interval with extremities in $Z$ contains points that are not in $Z$. Finally, preimages of branches are exactly intervals between points in $Z$ and the statement follows.
\end{proof}

\begin{remark}The dendrite $D_3$ is homeomorphic to the \emph{continuum self-similar tree} described in \cite{bonk2018continuum}. In this reference, it is proved that the \emph{continuous random tree} of Aldous is almost surely homeomorphic to $D_3$ as well.

The continuum self-similar tree $T$ is the attractor of the dynamical system of the plane generated by the following three maps

$$f_1(z)=\frac{1}{2}(z-1),\ f_2(z)=\frac{1}{2}(\overline{z} +1), f_3(z)=\frac{1}{2}(\overline{z}+i).$$

This means that $T$ is the unique compact subset of $\C$ such that $T=f_1(T)\cup f_2(T)\cup f_3(T)$. It is also obtained in the following way. One starts with $T_0=[-1,0]\cup[0,i]\cup[0,1]$ and defined $T_{n+1} =f_1(T_n)\cup f_2(T_n)\cup f_3(T_n)$. Finally, $T$ is the closure of the union of all the $T_n$’s.

\begin{figure}
\begin{center}
\begin{tikzpicture}[scale=4,ultra thin]
\draw(-1,0)--(1,0);
\draw(0,0) node[below] {$0$};
\draw(0,1) node[above] {$i$};
\draw(1,0) node[below] {$1$};
\draw(-1,0) node[below] {$-1$};
\input{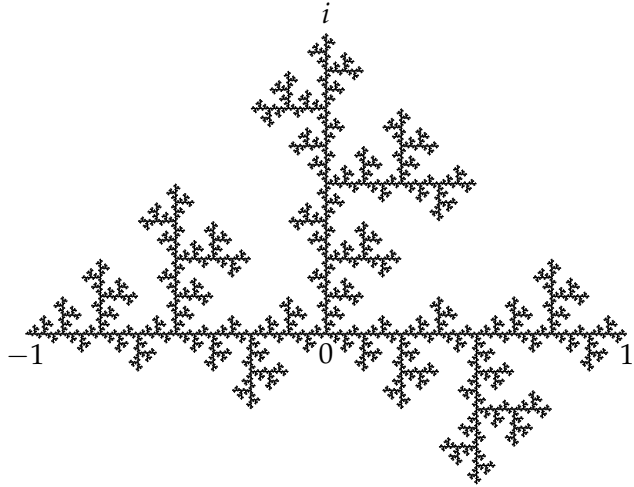}
\end{tikzpicture}
\caption{The continuum self-similar tree.}
\end{center}
\end{figure}\end{remark}
%
%
%

\subsection{Kaleidoscopic Colorings}

Let us recall the notion of kaleidoscopic coloring for dendrites introduced in \cite{MR4017014}. Let us denote by $\widehat{D_n}$ the set of branches of the dendrite $D_n$, that is the collection of all connected complements of some branch points and let $\widehat{x}$ be the collection of branches around the branch point $x$. Let us denote $[n]=\{1,...,n\}$. A \emph{coloring} is a map $c\colon \widehat{D_n}\to [n]$ such that the restriction of $c$ on branches around a given branch point is a bijection. Such a coloring is \emph{kaleidoscopic} if for any $x\neq y\in\Br(D_n)$ and $i\neq j\in[n]$ there is a branch point $z\in]x,y[$ such that $c(B_z(x))=i$ and $c(B_z(y))=j$. 
 \begin{center}
\setlength{\unitlength}{0.5cm}
\thicklines
\begin{picture}(4,4)
\cbezier(2,2)(2.25,3)(3,3.75)(4,4)
\cbezier(2,2)(2.25,1)(3,.25)(4,0)
\cbezier(2,2)(1.75,3)(1,3.75)(0,4)
\cbezier(2,2)(1.75,1)(1,.25)(0,0)
\put(2,2){\circle*{0.2}}\put(1.8,1){$z$}
\put(0,2){\circle*{0.2}}\put(-.2,1.5){$x$}
\put(4,2){\circle*{0.2}}\put(3.8,1.5){$y$}
\put(.5,2.8){$i$}
\put(3.3,2.8){$j$}
\end{picture}
\end{center}

We refer to \cite[\S3]{MR4017014} for details and in particular, existence (and abundance) of such kaleidoscopic colorings.

The group of homeomorphisms of the dendrite, $\Homeo(D_n)$, acts faithfully on the set of branch points of the dendrite and the image of the homomorphism to the symmetric group $\Sym(\Br(D_n))$ is actually a closed subgroup with respect to the pointwise convergence topology. This topology on $\Homeo(D_n)$ actually coincides with the uniform converge on the dendrite \cite[Proposition 2.4]{MR4017014}. 
Let $c$ be a kaleidoscopic coloring of $\widehat{D_n}$, $x$ be a branch point and $g\in\Homeo(D_n)$. The \emph{local action} of $g$ at $x$ is the element $\sigma_c(g,x)\in\Sym([n])$ such that the following diagram is commutative.

$$\begin{tikzcd}
\widehat{x} \arrow[r, "g"] \arrow[d, "c"]
& \widehat{g(x)} \arrow[d, "c"] \\ {[n]} \arrow[r, "{\sigma_c(g,x)}" ]
&{[n]} 
\end{tikzcd}$$
That is $\sigma_c(g,x)=c\circ g|_{\hat{x}}\circ c^{-1}$. For a subgroup $\Gamma\leq \Sym([n])$, the \emph{kaleidoscopic group} $\mathcal{K}(\Gamma)$ is the subgroup of elements $g\in\Homeo(D_n)$ such that $\sigma_c(g,x)\in\Gamma$ for all $x\in\Br(D_n)$. 
\section{Definition of the action on the circle and consequences}

\subsection{Definition of the action}Let us start with an observation on the set with three points $[3]=\{1,2,3\}$ and let us consider it with the cyclic order induced by the standard linear order on the natural numbers. The stabilizer of this cyclic order in the symmetric group $S_3$ is exactly $A_3$, the alternating group, which is the cyclic group $C_3$ as well. \\

In the remaining of this note, we identify $D_3$ with the Julia set of $z\mapsto z^2+i $ described in the previous section. The identification gives a cyclic order on the  three branches around each branch point of $D_3$ in the following way.

Let $x$ be a branch point of $D_3$, we denote by $z_1,z_2,z_3\in\bS^1$ its preimages positively ordered. Each of the three branches around $x$ has a preimage by the Carathéodory loop $\varphi$ that is a  positive interval $(z_i,z_{I+1})_+$. Let $A,B,C$ be the three branches around $x$. We define $i(A)$ to be the unique $i\in[3]$ such that $\varphi((z_i,z_{i+1})_+)=A$ and similarly, we define $i(B)$ and $i(C)$. We say that $(A,B,C)$ is positively ordered if $(i(A),i(B),i(C))$ is so.

More concretely, let us define $o_\varphi(A,B,C)=+$ if $(i(A),i(B),i(C))$ is positively ordered and $o_\varphi(A,B,C)=-$ otherwise. Similarly for a coloring $c$ of $\widehat{D_3}$, let us define $o_c(A,B,C)=+$ if $(c(A),c(B),c(C))$ is positively ordered in $[3]$ and $o_c(A,B,C)=-$ otherwise.

\begin{lemma} Let $x,y$ be distinct branch points in $D_3$ with respective preimages $x_1,x_2,x_3$ and $y_1,y_2,y_3$ by $\varphi$ and such that $(x_1,x_2,x_3,y_1,y_2,y_3)$ is cyclically ordered in $\bS^1$. There are $z_+$ and $z_-$ in $\Br(D_3)\cap(x,y)$ such that 
$$\varphi^{-1}(B_+)\subset[x_3,y_1]_+,$$
$$\varphi^{-1}(B_-)\subset[y_3,x_1]_+,$$
where $B_\pm$ is the branch at $z_\pm$ that does not contain $x$ nor $y$.
\end{lemma}

\begin{proof} By Lemma~\ref{branches}, the interval $(x_3,x_1)_+$ is the preimage of the branch at $x$ that contains $y$ and similarly, $(y_3,y_1)_+$ is the preimage of the branch at $y$ that contains $x$. Since $D(x,y)$ is the intersection of these two branches, $$\varphi^{-1}(D(x,y))=(x_3,x_1)_+\cap(y_3,y_1)_+=(x_3,y_1)_+\cup(y_3,x_1)_+.$$

Let $z$ be a branch point in $(x,y)$. Since the preimage by $\varphi$ of the branch at $z$ that does not contain $x$ nor $y$ is an interval in $(x_3,y_1)_+\cup(y_3,x_1)_+$, it is included in $(y_3,x_1)_+$ or in $(x_3,y_1)_+$. 

The unique critical point of $z\mapsto z^2+i$ is $z=0$. So, this map is a local diffeomorphism except at $z=0$, which is a regular point of the dendrite.  Since $\varphi$ semi-conjugates $z\mapsto z^2$ and $z\mapsto z^2+i$, the order of the points $\varphi(z)$ and $\varphi(z^2)$ are the same as soon as $\varphi(z)\neq0$. Actually, the images of $(y_3,x_1)_+$ under the iteration of $z\mapsto z^2$ cover the whole circle. So the intervals  $(y_3,x_1)_+$ and $(x_3,y_1)_+$ contain preimages of branches  at points $z\in(x,y)$ and we find points $z_\pm$ as in the statement.

\end{proof}

\begin{proposition}\label{coloring}There is a kaleidoscopic coloring $c$ of $\widehat{D}_3$ such that for any distinct branches around the same branch point $A,B,C$, distinct branches around the same branch point, $o_\varphi(A,B,C)=+$ if and only if  $o_c(A,B,C)=+$.\end{proposition}

\begin{proof} Such a coloring is constructed by induction. Let us first enumerate the points in $\Br(D_3)=\{b_n\}_{n\in\N}$.  We construct sequences of finite subsets $B_n\subset \Br(D_3)$ and maps $c_n\colon \widehat{B_n}\to[3]$ with the following properties for each $n\in\N$
\begin{enumerate}
\item $B_{n}\subset B_{n+1}$, $b_n\in B_{n}$,
\item the restriction of $c_{n+1}$ to $B_n$ coincides with $c_n$,
\item for each $b\in B_n$, $c_n|_{b}\colon \widehat{b}\to[3]$ is a bijection,
\item for all $b\in B_n$, $o_{c_n}$ and $o_{\varphi}$ coincide at $\widehat{b}$,
\item for all $b\neq b’$  in $B_{n-1}$ and $i\neq j\in\{1,2,3\}$, there is $b^{\prime \prime}\in B_n\cap(b,b’)$ such that $c_n(B_{b^{\prime\prime}}(b))=i$ and $c_n(B_{b^{\prime \prime}}(b’))=j$.
\end{enumerate}

Assume such sequences have been constructed. We simply define $c$ on $\widehat{b_n}$ to be $c_n|_{\widehat{b_n}}$. Thanks to Items (1) and (2), $c$ is well defined on $\widehat{D_3}$, it is a coloring thanks to Item (3), it is kaleidoscopic thanks to Item (5) and the cyclic orders $o_c$ and $o_\varphi$ coincide thanks to Item (4).\\

For a subset $Y\subset \Br(D_3)$, and a coloring $c\colon \widehat{Y}\to[3]$, we denote by $o_{c,b}$
 the cyclic order induced by $c$ at $b\in Y$ and similarly we denote by $o_{\varphi,b}$ the cyclic order induced by $\varphi$ at $b$. Let us define $B_0$ to be $\{b_0\}$ and define $c_0$ on $\widehat{b_0}$ arbitrarily such that $o_{c_0,b_0}$ and $o_{\varphi,b_0}$ coincide. Assume $B_n$ and $c_n$ have been constructed. Let us define $B’_n=B_n\cup\{b_{n+1}\}$ and $c’_n $ to coincide with $c_n$ on $\widehat{B_n}$ and to be arbitrarily defined on $\widehat{b_{n+1}}$ arbitrarily such that  $o_{c,b_{n+1}}$ and $o_{\varphi,b_{n+1}}$ coincide if  $b_{n+1}\notin B_n$.

Now for two points $x\neq y\in B_n$ and $i\neq j\in \{1,2,3\}$, let $k\in\{1,2,3\}$ distinct from $i$ and $j$. If $(i,k,j)$ is positive we choose $z_+\in (x,y)\cap\Br(D_3)$ given Lemma~\ref{branches}, we add $z_+$ to $B_{n+1}\supset B_n$ and define $c_{n+1}$ at $z_+$ by $c_{n+1}(B_{z_+}(x))=i$, $c_{n+1}(B_{z_+}(y))=y$ and $c_{n+1}(B)=k$ where $B$ is the third branch at $z_+$. One has $o_{c,z_+}(B_{z_+}(x),B,B_{z_+}(y))=+=o_{\varphi,z_+}(B_{z_+}(x),B,B_{z_+}(y))$. Similarly, if $(i,k,j)$ is negative, we choose $z_-\in (x,y)\cap\Br(D_3)$ given Lemma~\ref{branches}, we add $z_-$ to $B_{n+1}$ and define $c_{n+1}$ at $z_-$ by $c_{n+1}(B_{z_-}(x))=i$, $c_{n+1}(B_{z_-}(y))=y$ and $c_{n+1}(B)=k$ where $B$ is the third branch at $z_+$. By construction, $o_{c,z_\pm}$ and $o_{\varphi,z_\pm}$ coincide.\end{proof}

Let $G=\mathcal{K}(A_3)$ be the kaleidoscopic group associated to the coloring $c$ constructed in Proposition~\ref{coloring} and the group $A_3$. Let us now define an orientation-preserving action of $G$ on $\bS^1$. This action will be defined in two steps. First we define it on preimages of branch points via $\varphi$ and second, density of these points and order preserving will show that it extends uniquely continuously on the whole of $\bS^1$. The associated homomorphism $\rho\colon G\to\Homeo^+(\bS^1)$ will make the following diagram commutative, that is for all $z\in \bS^1$, $\varphi(\rho(g)z)=g\varphi(z)$.

$$\begin{tikzcd}
\bS^1 \arrow[r, "\rho(g)"] \arrow[d, "\varphi"]
& \bS^1 \arrow[d, "\varphi"] \\ D_3 \arrow[r, "g" ]
&D_3 \end{tikzcd}$$

Let us fix $g\in G$. We will define a bijection $\rho(g)$ on the set $\tilde B$ of preimages of branch points in $\bS^1$. Let $x$ be a branch point of $D_3$. Let $B_1,B_2,B_3$ be the branches around $x$ such that $c(B_i)=i$. Similarly, let $x’=gx$ and $B’_i$ the colored branches around $x’$ with $c(B_i')=i$. The local action $\sigma$ satisfies $B’_{\sigma(g,x)(i)}=g(B_i)$. If we set $z_1,z_2,z_3$ to be the preimages of $x$ such that $\varphi^{-1}(B_i)=(z_i,z_{i+1})_+$ (cyclically understood) and similarly for the preimages $z’_1,z’_2,z’_3$ of $x’$ then we define $\rho(g)(z_i)$ to be $z’_{\sigma(g,x)(i)}$. In particular, once $\rho(g)$ is entirely defined on $\tilde B$,  by Proposition~\ref{coloring}, 
\begin{equation}\label{cloclo}\rho(g)\left((z_i,z_{i+1})_+\cap\tilde B\right)=\left(z’_{\sigma(g,x)(i)},z’_{\sigma(g,x)(i+1)}\right)_+\cap \tilde B.\end{equation}

\begin{lemma} The map $\rho(g)\colon \varphi^{-1}\left(\Br(D_3)\right)\to \varphi^{-1}\left(\Br(D_3)\right)$ is a bijection such that  $\varphi(\rho(g)z)=g\varphi(z)$ for all $z\in \varphi^{-1}\left(\Br(D_3)\right)$ and which preserves the cyclic order on the circle.\end{lemma}

\begin{proof} The first point is straightforward from the definition of $\rho(g)$. Let us prove that $\rho(g)$ preserves the cyclic order on $ \varphi^{-1}\left(\Br(D_3)\right)$. Let $z_1,z_2,z_3\in  \varphi^{-1}\left(\Br(D_3)\right)$ positively ordered. If these points have the same image under $\varphi$ then their images $\rho(g)(z_i)$ are positively ordered as well by construction. 

Let us assume their images are two distinct points $x_1,x_2\in D_3$.  Up to reorder cyclically the $z_i$’s, we may assume that $\varphi(z_1)=\varphi(z_2)=x_1$ and $\varphi(z_3)=x_2$. Let $B_1$ be the branch around $x_1$ that contains $x_2$ and let $B_2$, $B_3$ the two other branches ordered cyclically with above cyclic order. Let $y_1,y_2,y_3$ be the preimages of $x_1$ by $\varphi$ such that $\varphi^{-1}(B_{i+1})=(y_i,y_{i+1})_+$. Since $z_3\in(y_3,y_1)_+$, $z_1=y_1$ and $z_2=y_2$ or $y_3$, or $z_1=y_2$ and $z_2=y_3$. Let $B’_1$ be the branch around $g(x_1)$ that contains $g(x_2)$ and $B_2’, B_3’$ the two other branches cyclically ordered. Since $g(B_1)=B_1'$ and $\sigma_c(g,x_1)\in A_3$, $g(B_2)=B’_2$ and $g(B_3)=B_3’$. Let  $y’_1,y’_2,y’_3$ be the preimages of $g(x_1)$ such that $\varphi^{-1}(B_{i+1})=(y’_i,y’_{i+1})_+$. By construction, $\rho(g)(y_i)=y’_i$ and by Equation \eqref{cloclo}, $(\rho(g)(z_1),\rho(g)(z_2),\rho(g)(z_3))$ is positively ordered.

Let us assume that the images of $z_1,z_2,z_3$ are three distinct points: $x_1,x_2,x_3$. There are two possibilities: these points are aligned or they are not and this last case, they have a center $c$ which is different from the $x_i$’s. In this last case, we can use the three preimages $z’_1,z_2’,z_3’$ of $c$ labeled such that  $\varphi^{-1}(B_i)=(z_i’,z_{i+1}’)_+$ where $B_i$ is the branch around $c$ containing $x_i$. By the first point of this proof, the images of  $z’_1,z_2’,z_3’$ by $\rho(g)$ are positively ordered and thus the same hold for the images of $z_1,z_2$ and $z_3$.

Let us conclude with the case where $x_1,x_2,x_3$ are distinct and aligned. Up to permute cyclically these points, we may assume that $x_2\in[x_1,x_3]$. Let $z’_3$ and $z’_1$ the preimages of $x_2$ distinct from $z_2$ such that $(z_2,z’_3,z’_1)$ is positively ordered in the circle. Considering the three possible relative positions of $z_1, z_2$ with respect to the positive intervals with extremities $z’_1,z_2,z’_3$ and thanks to Equation \eqref{cloclo}, we see that the images of $z_1,z_2,z_3$ are positively ordered as well.
\end{proof}

Since $\varphi^{-1}\left(\Br{D_3}\right)$ is dense in $\bS^1$ and the topology on the circle is defined by this cyclic order, we get the following.

\begin{proposition} The map $\rho(g)$ extends uniquely to an element of $\Homeo^+(\bS^1)$ and it satisfies  $\varphi(\rho(g)z)=g(\varphi(z))$ for all $z\in \bS^1$.\end{proposition}

The proof of this proposition is a direct consequence of the following easy lemma that we will use later again.

\begin{lemma}\label{cyclic_homeo} Let $D\subset \bS^1$ be a dense subset and $f\colon D\to \bS^1$ be an injective map that preserves the cyclic order and with dense image.  Then $f$ extends uniquely to a homeomorphism of $\bS^1$.
\end{lemma}

\begin{proof}Since $f$ preserves the order, if $z_1,\dots, z_4\in D$ and $[z_2,z_3]_+\subset [z_1,z_4]_+$ then $f(D\cap[z_2,z_3]_+)\subset [f(z_1),f(z_4)]_+$. Since $D$ is dense, for any point $z\in\bS^1$,  $$\{z\}=\bigcap_{z_-,z_+\in D,\ z\in[z_-,z_+]_+}[z_-,z_+]_+.$$
By density of $f(D)$, the compact intersection $\displaystyle \bigcap_{z_-,z_+\in D,\ z\in[z_-,z_+]_+}[f(z_-),f(z_+)]_+$ is reduced to a point. We extend $f$ to $\bS^1$ by defining $f(z)$ to be this point. This extension preserves the cyclic order on $\bS^1$. Actually if $z\in[z_-,z_+]_+$ and these points are distinct then one can find $z_1,\, z_2,\, z_3\in D$ such that $z_-\in(z_1,z_2)_+, z\in(z_2,z_3)_+$ and $z_+\in (z_3,z_1)_+$. By construction  $f(z_-)\in(f(z_1),f(z_2))_+, f(z)\in(f(z_2),f(z_3))_+$ and $f(z_+)\in (f(z_3),f(z_1))_+$ and thus $(f(z_-),f(z),f(z_+))$ is positively ordered. This extension is injective since for any $z\neq z'$, one can find $(z_1,z_2,z_3)\in D^3$ positively ordered such that $z\in(z_1,z_2)_+$ and $z'\in(z_2,z_3)_+$. Thus their images have the same orders and so $f(z)\neq f(z')$. Moreover $f$ is continuous since it preserves the cyclic order and it is surjective by compactness and density of its image. The density of $D$ implies the uniqueness of this extension.
\end{proof}
\subsection{Proof of the statements}

\begin{proposition} The map $\rho\colon  G\to\Homeo^+(\bS^1)$ is a topological group embedding.\end{proposition} 

\begin{proof} The map $\rho$ is a group homomorphism because $\sigma_c$ is a cocycle. Let $(g_n)$ be a sequence of elements of $G$ converging to the identity in $G$ then $g_n(x)=x$ for all $x\in \Br\left(D_3\right)$ and $n$ large enough so the same hold for $\rho(g_n)(z)$ and $z\in \varphi^{-1}\left(\Br{D_3}\right)$ because by choosing points $x'$ in branches around $x$, we see that $\sigma_c(g_n,x)=\Id$ for $n$ large enough. By compactness of $\bS^1$, this implies that $\rho(g_n)$ converges uniformly to the identity. So $\rho$ is a continuous homomorphism. Conversely, if $\rho(g_n)$ converges uniformly to the identity on $\bS^1$ then  the commutative diagram and the continuity of $\varphi$ imply that $g_n$ converges uniformly to the identity in $G$. Since $G$ is Polish, this implies that $G$ is closed in $\Homeo^+(\bS^1)$.\end{proof}

The desired action $G\action \bS^1$ is now constructed. So, we can prove the existence of exactly 3 orbits and oligomorphy of this action.

\begin{proof}[Proof of Proposition \ref{orbits}] We claim that two points of $\bS^1$ are in the same orbit if and only if their images in the dendrite have same order. This is, of course, a necessary condition. In the other direction, this is a consequence of \cite[Lemma 5.8]{MR4017014}. 
\end{proof}

\begin{proof}[Proof of Proposition \ref{oligomorphy}]The oligomorphy follows from the oligomorphy of the action $G\action D_3$ \cite[\S 5D]{MR4017014} and the fact there are only finitely many subsets $S\subset\bS^1$ of size $n$ with the same image by $\varphi$. 
\end{proof}

To get non-elementarity, we prove something stronger, namely strong proximality. Let us recall that a continuous action of a group $\Gamma$ on a compact space $X$ is \emph{strongly proximal} if any probability measure on $X$ contains a Dirac mass in the weak*-closure of its $G$-orbit. 

\begin{lemma} The action $G\action \bS^1$ is minimal.\end{lemma}

\begin{proof} Let $K$ be some closed $G$-invariant non-empty subspace of $\bS^1$. By minimality of the action $G\action D_3$, $\varphi(K)=D_3$. So $K$ contains all preimages of end points and since these points are dense, $K=\bS^1$.\end{proof}

\begin{proposition}\label{proximality}The action $G \curvearrowright \bS^1$ is strongly proximal.\end{proposition}

\begin{proof} This follows from the fact that any interval of $\bS^1$ (with non-empty complement) can be sent to the interior of some other non-trivial interval. This follows from the fact that the image of an interval $[z_1,z_2]_+$ by $\varphi$ is contained in a branch around a branch point $x\notin \varphi([z_1,z_2]_+)$ can be sent to any given branch by an element of $G$. See \cite[Theorem 10.1]{DuchesneMonod} .\end{proof}

\begin{remark} The non-elementarity of the action also follows directly from \cite[Proposition 3.2]{DuchesneMonod}. Actually, if there is an invariant measure on $\bS^1$, its image by $\varphi$ is an invariant measure on $D_3$ and thus there is a finite orbit on $D_3$. This is not the case.\end{remark}

\begin{remark} Since end points have a unique antecedent, $\varphi$ induces a bijection between end points of $D_3$ and a subset, which is actually comeager, of $\bS^1$. This yields a cyclic order on the set of end points and allows to construct a bounded cohomology class on $G$, the so-called \emph{orientation class}. This is reminiscent of the construction of cocycles done in \cite[Proposition 9.1]{MR4017014}. \end{remark}

Let us prove Property (T) for $G$.
\begin{proof}[Proof of Proposition~\ref{T}] It is oligomorphic by  \cite[Theorem 5.9]{MR4017014} and thus the result follows from \cite[Theorem 1.1]{MR3417738}.\end{proof}
Let us conclude with the proof that $G$ is exactly the automorphism group of the laminational equivalence relation $\sim_\varphi$ induced by the Carathéodory loop.

\begin{proof}[Proof of Proposition~\ref{lamination}]
By definition of the action of $G$ on the circle, it is clear that the equivalence relation is invariant under this action. Let $g\in\Homeo^+(\bS^1)$ preserving the equivalence relation. Since $\varphi(\bS^1)=D_3$, this defines a map $\psi(g)\colon D_3\to D_3$ such that for all $z\in\bS^1$, $\psi(g)(\varphi(z))=\varphi(g(z))$. Counting antecedents of points in $D_3$, one sees that $\psi(g)$ maps branch points to branch points and preserve the betweenness relation  because if for three points $x,y,z\in D_3$, $x\in[y,z]\iff x\in\varphi\left([y',z']_+\right) \cap\varphi\left([z',y']_+\right)$. So by \cite[Proposition 2.4]{MR4017014}, $\psi(g)$ is a homeomorphism of the dendrite $D_3$. It preserves the cyclic order on branches around a branch point (because $g$ preserves the cyclic order of the circle)
 and thus it is an element of $G=\mathcal{K}(A_3)$. Moreover $\rho(\psi(g))=g$ and thus $\rho$ is an isomorphism between $G$ and the automorphism group of the laminational equivalence relation $\sim_\varphi$. \end{proof}

\section{Uniqueness of the action}

For an end point $\xi\in D_3$, we denote by $G_\xi$ its stabilizer, which is a closed subgroup. A topological group $H$ is \emph{extremely amenable} if any $H$-flow has a fixed point. This property may seem a priori very surprising since any locally compact group has a free action on a compact space (its Stone-\v{C}ech compactification) and thus no non-trivial locally compact group is extremely amenable. Nonetheless, this property is not so rare for very large group. For example, the unitary group of a separable Hilbert with its strong topology or the automorphism of the rational numbers with its canonical order are extremely amenable \cite{MR708367,MR1608494}.

For automorphism groups of countable structures, Kechris, Pestov and Todorcevic characterized when a Fraïssé limit has an extremely amenable automorphism group. It is the case if and only if this Fraïssé limit is a Fraïssé order class with the Ramsey property \cite{MR2140630}.  They used this result to identify some universal minimal flows of some automorphism groups of countable structures. We refer to their paper for an introduction to this subject and a precise meaning.

The following proposition is simply \cite[Corollary 7.3(i)]{MR2140630} once $G_\xi$ has been identified with the automorphism group of Fraïssé limit $\Aut(\mathbb{DT}_2)$ as it will be explained in the proof.

\begin{proposition} The topological group $G_\xi$ is extremely amenable.
\end{proposition}

\begin{proof} Following the strategy developed in \cite{MR2140630}, we prove that $G_\xi$ is the automorphism group of an order Fraïssé limit with the Ramsey property. The Ramsey property will come from \cite{MR3436366} and we try to follow notations there. The universe of this limit is the set of branch points of $D_3$. We endow it with the following partial order $a\leq b$ if $a\in[\xi,b]$. The cyclic order on $\bS^1$ induces a linear order on $\bS^1\setminus\{\varphi^{-1}(\xi)\}$ in the following way $z<w$ if $(\varphi^{-1}(\xi),z,w)$ is positively ordered. Now for $x\in \Br(D_3)$, we define $x_1,x_2,x_3\in\bS^1$ its three preimages by $\varphi$ such that $x_1<x_2<x_3$. We write $R_1(a,b)$ if $b\in\varphi((a_1,a_2)_+)$ and $R_2(a,b)$ if $b\in\varphi((a_2,a_3)_+)$. Thus, $R_1$ and $R_2$ are binary relations. These relations indicate in which branch at $a$, a successor $b$ for $\leq$ of $a$ lies. The remaining of the proof is devoted to show that $(\Br(D_3),\leq,R_1,R_2)$ is isomorphic to the Fraïssé limit $\mathbb{DT}_2$ which has the Ramsey property thanks to \cite[Theorem 6.1]{MR3436366} and that $G_\xi$ is the automorphism groups of $(\Br(D_3),\leq,R_1,R_2)$. 

Let us prove first this second point. It is clear that $G_\xi$ acts faithfully by automorphisms on $(\Br(D_3),\leq,R_1,R_2)$. Conversely, let $g\in \Aut(\Br(D_3),\leq,R_1,R_2)$. To show that $g$ is the restriction of some element $\tilde{g}$ on $\Homeo\left(\Br(D_3)\right)$, it suffices to prove that $g$ preserves the betweenness relation \cite[Proposition 2.4]{MR4017014}. Let $a,b,c\in\Br(D_3)$ such that $c\in [a,b]$. If $a\leq b$ then $a\leq c\leq b$. So $g(a)\leq g(c)\leq g(b)$ and thus $g(c)\in[g(a),g(b)]$. Let $a\circ b$ be the meet of $a$ and $b$ for $\leq$ that is $a\circ b=\inf\{a,b\}$ which the center between $\xi,a$ and $b$ in $\Br(D_3)$. Now, $c$ belongs to $[a\circ b,a]$ or $[a\circ b,b]$ and we repeat the above with $(a\circ b, a)$ or $(a\circ b,b)$ and thus  $g(c)\in[g(a),g(b)]$. So $g$ is the restriction of some $\tilde{g} \in \Homeo\left(\Br(D_3)\right)$. Let $(x_n)$ be a minimizing sequence in $\Br(D_3)$ for $\leq$. It converges to $\xi$ for the topology on  $D_3$. Its image by $\tilde{g}$ has the same property and thus $\tilde{g}(\xi)=\xi$.  

For $x\in \Br(D_3)$ with preimages $x_1,x_2,x_3$ ordered as before. The three branches around $x$ are $\varphi((x_3,x_1)_+)$ (the one containing $\xi$), $\varphi((x_1,x_2)_+)$ and $\varphi((x_2,x_3)_+)$. By construction, $\varphi((x_1,x_2)_+)=\left\{y\in \Br(D_3), y\geq x, R_1(x,y)\right\}$ and $\varphi((x_2,x_3)_+)=\{y\in \Br(D_3), y\geq x, R_2(x,y)\}$. Since the same fact holds mutatis mutandis for $g(x)$ as well, one as that $\sigma_{c}(\tilde{g},x)\in A_3$. Thus $\tilde{g}\in G_\xi$ and $g\mapsto \tilde{g}$ is an isomorphism between $\Aut(\Br(D_3),\leq,R_1,R_2)$ and $G_\xi$.

Now, let us prove that $(\Br(D_3),\leq,R_1,R_2)$ is isomorphic to $\mathbb{DT}_2$ which appears in \cite{{MR3436366}}. The former countable. structure is the Fraïssé limit of the class of finite structures $\mathcal{DT}_2$ where $(A,\leq^A,R^A_1,R^A_2)\in\mathcal{DT}_2$ if it satisfies the following properties;
\begin{enumerate}
\item $(A, \leq_A)$ is a meet semilattice (any pair of points $\{a,b\}$ has an infimum $a\circ b$);
\item for any $a\in A$, $\{b\in A, b\leq^Aa\}$ is linearly ordered with $\leq_A$, i.e., $(A, \leq^A)$ is a treeable semilattice.
\item Each $a\in A$ has at most 2 immediate successors for $\leq^A$.
\item $R^A_i(a,b)\implies a<^Ab$.
\item $a<^Ab\implies R_1(a,b)$ or $R_2(a,b)$.
\item $R^A_i(a,b)$ and $b\circ c=a\implies \neg R^A_i(a,c)$.
\end{enumerate}

By the characterization of Fraïssé limits, to prove that $(\Br(D_3),\leq,R_1,R_2)$ is isomorphic to $\mathbb{DT}_2$, it suffices to prove that any finite subsets of $\Br(D_3)$ is an element of $\mathcal{DT}_2$ and conversely any element of $\mathcal{DT}_2$ embeds in $(\Br(D_3),\leq,R_1,R_2)$. 

One check easily that any finite substructure of  $(\Br(D_3),\leq,R_1,R_2)$ belongs to $\mathcal{DT}_2$. Now, we prove by induction that any $(A,\leq^A,R_1^A,R^A_2)\in\mathcal{DT}_2$ embeds in $(\Br(D_3),\leq,R_1,R_2)$. We proceed by induction on the cardinal of $A$. Let $a\in A$ be the root (i.e., the minimal element in $(A,\leq^A)$), we choose any point $i(a)\in \Br(D_3)$. We set $A_1=\{b\in A,\ R_1^A(a,b)\}$ $A_2=\{b\in A,\ R_1^A(a,b)\}$. The substructures $A_1$ and $A_2$ are in $\mathcal{DT}_2$ with cardinality less than the one of $A$. So, we have an embedding  $i_1$ of $A_1$ in the first branch around $I(a)$ that does not contain $\xi$, which is isomorphic to  $(\Br(D_3),\leq,R_1,R_2)$ actually, and an embedding  $i_2$ of $A_2$ in the first branch around $I(a)$ that does not contain $\xi$. The concatenation of $i,i_1,i_2$ yields the embedding in $(\Br(D_3),\leq,R_1,R_2)$.  \end{proof}

\begin{proof}[Proof of Theorem \ref{unique}]
Let us assume that $G$ acts minimally on the circle and let us prove first that there is no invariant probability measure. To avoid ambiguity, let us call this action $\beta\colon G\times \bS^1\to\bS^1$ and let $\alpha$ be the action described in Theorem~\ref{action}. If there is an invariant probability measure $\mu$ then the rotation number is a homomorphism $G\to \R/\Z$ (see \cite[2.2.2]{MR2809110}. Since $G$ is simple, this homomorphism is trivial and the support of $\mu$ is a non-empty closed invariant subset consisting of fixed points. This is a contradiction with minimality of the action.

By extreme amenability of $G_\xi$, this group has a fixed point $x\in\bS^1$ for $\beta$. This allows us to define the orbit map $\psi\colon G/G_\xi\to\bS^1$ such that $\psi(gG_\xi)=\beta(g,x)$. Since $G/G_\xi\simeq \Ends(D_3)$, by an abuse of notation, we also write $\psi(\eta)=\psi(gG_\xi)$ is $\eta=g\xi$. Since $G$ acts doubly transitively on the end points of $D_3$, the group $G_\eta$ is maximal in $G$ for any $\eta\in\Ends(D_3)$. If $\psi(\eta)=\psi(\eta')$ then $G_\eta\cup G_{\eta'}\leq \Stab(\psi(\eta))$ and since $\Stab(\psi(\eta))<G$ then $\eta=\eta'$. So $\psi$ is injective.

The cyclic order on $\bS^1$ induces a cyclic order on $\Ends(D_3)$ via the Carathéodory loop $\varphi\colon \bS^1\to D_3$. More precisely, we define, for $\xi_1,\xi_2,\xi_3$ distinct, $\co(\xi_1,\xi_2,\xi_3)\in\{\pm\}$ according to the orientation of $(z_1,z_2,z_3)$ in the circle such that $\varphi(z_i)=\xi_i$.

We observe that for fixed $\xi\neq\eta\in\Ends(D_3)$ $\co(\xi,\lambda,\eta)=\co(\xi,\lambda',\eta)$ if and only if there is $g\in G$ fixing $\eta$ and $\xi$ such that $g(\lambda)=\lambda'$. Actually, there are two orbits for the action of $G$ on triple of distinct points in $\Ends(D_3)$. Two such triples are in the same orbit if and if they have the same cyclic orientation. More precisely, let $(\xi,\lambda,\eta)$ and $(\xi',\lambda',\eta')$ be two positively oriented triples in $\Ends(D_3)$. Let $c,c'$ be their respective centers. By \cite[Corollary 4.5]{MR4017014}, we found a homeomorphism $f_\xi\colon \overline{D(\xi,c)}\to\overline{D(\xi',c')}$ such that $f_\xi(\xi)=\xi'$ and $f_\xi(c)=c'$ preserving the coloring. Similarly, we construct $f_\lambda$ and $f_\eta$. Thanks to the patchwork lemma, \cite[Lemma 2.9]{MR3894039} we patch $f_\xi,f_\lambda$ and $f_\eta$ together to get an element $f\in G$ such that $f(\xi)=\xi',\ f(\lambda)=\lambda'$ and $f(\eta)=\eta'$.

Since $G$ is a simple group, the image of $G$ in $\Homeo(\bS^1)$ induced by the action $\beta$, is actually in $\Homeo^+(\bS^1)$.

Choose $(\eta,\lambda,\xi)\in\Ends(\bS^1)^3$ such that $\co(\eta,\lambda,\xi)=+$ and fix a cyclic orientation $\co_{\bS^1}$ on $\bS^1$ such that   $\co_{\bS^1}(\psi(\eta),\psi(\lambda),\psi(\xi))=+$. Observe that $G$ preserves $\co_{\bS^1}$ as well. So by transitivity of the action of $G$ of triple in $\Ends(D_3)$ with a given orientation, one has that $\co_{\bS^1}(\psi(\eta),\psi(\lambda),\psi(\kappa))=\co(\eta,\lambda,\kappa)$ for any triple of distinct points $(\eta,\lambda,\kappa)\in\Ends(D_3)$. 

In particular the map $\psi\circ\varphi\colon \varphi^{-1}(\Ends(D_3)\to\bS^1$ is an orientation-preserving map between two dense subsets of $\bS^1$. Thus it extends to a unique continuous $G$-map $h\colon\bS^1\to\bS^1$, which a homeomorphism by Lemma~\ref{cyclic_homeo}. This homeomorphism conjugates the action $\alpha$ and $\beta$.
\end{proof}

\begin{proof}[Proof of Theorem \ref{no smooth}] Let $\omega$ be a modulus of continuity and  let $\alpha$ be any minimal continuous action of $G$ by homeomorphisms on $\bS^1$. This action $\alpha$ is conjugated to the action described in Theorem~\ref{action}  by some homeomorphism $h\colon\bS^1\to\bS^1$. Let $\omega_h$ and $\omega_{h^{-1}}$ be the modulus of continuity of $h$ and $h^{-1}$. 

As before, we denote by $\rho\colon G\to\Homeo^+(\bS^1)$ the homomorphism by the action $G\action \bS^1$ described in Theorem~\ref{action} and we denote by $\rho_\alpha\colon G\to\Homeo^+(\bS^1)$, the one coming from the action $\alpha$, that is $\rho(g)=h\rho_\alpha(g)h^{-1}$.

For $g\in G$, if $\rho_\alpha(g)$ admits $\lambda \omega$ as modulus of continuity with $\lambda>0$ then $\rho(g)$ admits $\omega_{\lambda,h}$ as modulus of continuity where  $\omega_{\lambda,h}(r)=\omega_h(\lambda\omega(\omega_{h^{-1}}(r)))$.  So it suffices to prove that there is $g\in G$ such there is $x\in\bS^1$ and that for any $n\in \N$, there is $x_n$ such that \begin{equation}\label{nomc}
d(\rho(g)(x),\rho(g)(x_n))>\omega_{n,h}(d(x,x_n)).
\end{equation}
 Let us construct such $g\in G$ by patch working (See \cite[Lemma 2.9]{MR3894039}). We fix two branch points $b,b'\in D_3$. We define $g$ to be trivial outside $D(b,b')$. Let $b_1,b_2,b_3$ and $b'_1,b'_2,b'_3$ the antecedents of $b$ and $b'$ cyclically ordered. Choose $x'_1\in(b'_3,b_1)_+$ whose image in the dendrite is an end point.  Let $c'_1\in D_3$ be the center of the points $b,\varphi(x'_1),b'$. Now, by definition of the kaleidoscopic coloring, we may find a point $x_1\in(x'_1,b_1)_+$ in the circle such that
\begin{itemize}
\item $\varphi(x_1)$ is an end point,
\item the center $c_1$ of $b,b'$ and $\varphi(x_1)$ satisfies $c\left(B_{c_1}(b)\right)=c\left(B_{c'_1}(b)\right)$ and $c\left(B_{c_1}(b')\right)=c\left(B_{c'_1}(b')\right)$,
\item $d(b_1,x'_1)>\omega_{n,h}(d(b_1,x_1))$.
\end{itemize}

Thanks to Corollary 4.5 in \cite{MR4017014}, there are homeomorphisms $f_1\colon \overline{D(c_1,b')}\to \overline{D(c_1',b')}$ and $h_1\colon \overline{D(\varphi(x_1),c_1)}\to \overline{D(\varphi(x'_1),c'_1)}$ that preserve the kaleidoscopic coloring. We define $g$ to coincide with $f_1$ and $h_1$ on $\overline{D(c_1,b')}$ and $\overline{D(\varphi(x_1),c_1)}$. 

Assume that the points $c_k,c'_k, x_k,x'_k,$ have been constructed and $g$ has been defined on $\overline{D(c_k,b')}$ and $\overline{D(\varphi(x_k),c_k)}$ for $k\leq n$. Choose $x'_{n+1}\in[x'_{n},b_1]$ whose image via the Carathéodory loop $\varphi$ in the dendrite is an end point in the dendrite.  Let $c'_{n+1}\in D_3$ be the center of the points $b,\varphi(x'_{n+1}),b'$. Now, by definition of the kaleidoscopic coloring, we may find a point $x_{n}\in(x'_{n+1},b_1)_+$ in the circle such that
\begin{itemize}
\item $\varphi(x_{n+1})$ is an end point,
\item $d(x_{n+1},b)<\frac{1}{n+1}$,
\item the center $c_{n+1}$ of $(b,b',\varphi(x_{n+1}))$ satisfies $c\left(B_{c_{n+1}}(b)\right)=c\left(B_{c'_{n+1}}(b)\right)$ and $c\left(B_{c_{n+1}}(b')\right)=c\left(B_{c'_{n+1}}(b')\right)$,
\item $d(b_1,x'_{n+1})>\omega_{n+1,h}(d(b_1,x_{n+1}))$.
\end{itemize}
There are homeomorphisms $f_{n+1}\colon \overline{D(c_{n+1},c_n)}\to \overline{D(c_{n+1}',c_n')}$ and $h_{n+1}\colon \overline{D(\varphi(x_{n+1}),c_{n+1})}\to \overline{D(\varphi(x'_{n+1}),c'_{n+1})}$ that preserve the kaleidoscopic coloring. We define $g$ to coincide with $f_{n+1}$ and $h_{n+1}$ on $ \overline{D(c_{n+1},c_n)}$ and $\overline{D(\varphi(x_{n+1}),c_{n+1})}$.  Since $x_n\to b$ and thus $c_n\to b$, $g$ is defined everywhere and the patchwork lemma \cite[Lemma 2.9]{MR3894039} proves that $g$ is a well-defined homeomorphism of $D_3$ and by construction, it preserves the kaleidoscopic coloring, that is, $g\in\mathcal{K}(\{Id\})<G$.  Moreover, by construction $\rho(g)(x_n)=x'_n$. In particular, if we set $x=b_1$ then $\rho(g)$ satisfies Equation~\ref{nomc} for all $n$ and the theorem is proved.
\end{proof}

\section{Universal minimal flow}
We endow the quotient space $G/G_\xi$ with the uniform structure induced by the left uniform on $G$. A basis of entourages is given by 
$$U_V=\{(gG_\xi,hG_\xi),\ h\in vgG_\xi, \textrm{with}\ v\in V\}$$ where $V<G$ is the pointwise stabilizer of a finite set of branch points. We denote by $\widehat{G/G_\xi}$ the completion of this uniform structure.

For the homeomorphism group of  generalized Wa\.zewski dendrites, identifications of the universal minimal flows have been done in \cite{MR3893292,duchesne2019topological}. When there is a unique type of branch points, say of order $n$, this corresponds to the kaleidoscopic group $\mathcal{K}(S_n)$. For the kaleidoscopic group $\mathcal{K}(A_3)$, the situation is a bit different but rely on similar ideas.

\begin{proof}[Proof of Theorem~\ref{umf}]  We claim that the space $\widehat{G/G_\xi}$ is compact or equivalently $G_\xi$ is  coprecompact. This means that for any open neighborhood of the identity $V\subset G$ there is a finite set $F\subset G$ such that $VFG_\xi=G$.

A basis of open neighborhoods  of the identity in $G$ is given by pointwise stabilizers of a finite set of branch points. Let $B$ such finite set of branch points, we may assume that $B$ center-closed, that is for any triple points in $B$, their center is still in $B$. A connected component $C$ of $D_3\setminus B$ is of two types: either its boundary is a point or it has two points. If there were more than 2 points, the center of three of them is in $C$ and we have a contradiction with the fact that $B$ is center-closed. There are finitely many such connected components because any element in $B$ has finite order (3 actually). 

For each such connected component $C$, we choose elements $f\in G$ in the following way. If $C$ has a unique point in its boundary, we choose $f\in G$ such that $f(\xi)\in C$. If $C$ has two points $a,b$ in its boundary, we choose elements $f_-,f_+\in G$ such that $f_-(\xi),\, f_+(\xi)\in C$, $(c(B_z(a)),c(B_z(f_-(\xi)),c(B_z(b)))$ is a negative triple in $\{1,2,3\}$ where $z$ is the center of $(a,b,f_-(\xi))$ and $(c(B_y(a)),c(B_y(f_+(\xi)),c(B_y(b)))$ is a positive triple in $\{1,2,3\}$ where $y$ is the center of $(a,b,f_+(\xi))$. We define $F$ to the finite set  of all these elements $f$ or $f_-,f_+$ for all connected component of $D_3\setminus B$.

Now let $g\in G$, $g(\xi)$ belongs to some connected component $C$ of $D_3\setminus B$. If $C$ has a unique point $a$ in its boundary, let $f$ be the unique element in $F$ such that $f(\xi)\in C$. By  Corollary 4.5 in \cite{MR4017014}, there is a homeomorphism $v_0$ of $C\cup\{a\}$ fixing $a$, such that $v_0(f(\xi))=g(\xi)$ and $v_0$ preserves the coloring, we extend $v_0$ to an element $v\in G$ that is trivial outside $C\cup\{a\}$. By construction $v\in V$ and $vf(\xi)=g(\xi)$.

If $C$ has two points $a,b$ in its boundary, we choose $f\in F$ to be the element such that $f(\xi)\in C$ and $(c(B_z(a)),c(B_z(g(\xi)),c(B_z(b)))$ and $(c(B_y(a)),c(B_y(f(\xi)),c(B_y(b)))$ have same cyclic order where $z$ is the center of $(a,g(\xi),b)$ and $y$ the one of $(a,f(\xi)),b)$. Now, by applying Corollary 4.5 in \cite{MR4017014} three times, we find $v_1\colon \overline{D(a,z)}\to\overline{D(a,y)}$, $v_2\colon \overline{D(b,z)}\to\overline{D(b,y)}$ and $v_3\colon \overline{B_z(g(\xi))}\to\overline{B_y(f(\xi))}$ that preserves the coloring and maps respectively $(a,z)$ to $(a,y)$, $(b,z)$ to $(b,y)$  and $(z,g(\xi))$ to $(y,f(\xi))$. We patch these three elements to get a homeomorphism $v_0$ of $\overline{D(a,b)}$. Observe that the local action $\sigma_c(v_0,x)$ is trivial  for any $x\in\Br(D_3)\cap D(a,b)$ except maybe for $x=z$ where  $\sigma_c(v_0,z)\in A_3$. We extend $v_0$ trivially outside $\overline{D(a,b)}$ to get an element $v\in V$ such that $vf(\xi)=g(\xi)$.

In both cases, this proves that $g\in VFG_\xi$ and the coprecompactness is proved.\\

Since $G_\xi$ is extremely amenable, for any $G$-flow, $G_\xi$ has a fixed point $x$ and the orbit map $g\mapsto gx$ induces a uniformly continuous map $G/G_\xi\to X$ that extends to a $G$-map $\widehat{G/G_\xi}\to X$.  It remains to prove that $\widehat{G/G_\xi}$ is minimal. We rely on the characterization of minimality obtained in \cite[Proposition 6.6]{zucker_2020}. It is equivalent to prove that $G_\xi$ is \emph{pre-syndetic} in $G$ (as coined in \cite[Definition 8.1]{basso2020topological}). This means that for any neighborhood of the identity $V$ there is a finite set $F$ such that $G=FVG_\xi$. 

So let us consider a finite set of branch points $B$ and let $V$ be its pointwise stabilizer. Let $b$ be the projection of $\xi$ on the smallest subdendrite containing $B$.  Let $C_1$ be the branch around $b$ containing $\xi$. Let us order cyclically $C_2,C_3$ the other branches around $b$ and choose $f_1,f_2,f_3\in G$ fixing $b$ such that $f_i(C_1)=C_i$ (using \cite[Corollary 4.5]{MR4017014} as before).  Let $F=\{f_1,f_2,f_3\}$.

Let $g\in G$.  Take $i\in\{1,2,3\}$ such that $g(\xi)\in C_i$. Thus $f_i^{-1}(g(\xi))\in C_1$. There is $v\in G$ trivial outside $C_1$ (and thus $v\in V$) such that $v(\xi)=f_i^{-1}(g(\xi))$ (as an application of \cite[Corollary 4.5]{MR4017014} again). So $v^{-1}f_i^{-1}g\in G_\xi$, that is $g=f_ivh$ for some $h\in G_\xi$. Thus $G=FVG_\xi$ and the minimality of $\widehat{G/G_\xi}$ is proved.
\end{proof}

\begin{corollary} Any $G$-flow is strongly proximal.\end{corollary}

\begin{proof} It suffices to prove that $\widehat{G/G_\xi}$ is strongly proximal. Actually, Theorem 10.1 in \cite{DuchesneMonod} shows that for any probability measure $\mu$ on $D_3$, one can find $g_n\in G$ such that $(g_n)_*\mu\to\delta_\xi$. Let $\nu$ be some probability measure on $\widehat{G/G_\xi}$ and let us apply the aforementioned result to $\mu=\pi_*(\nu)$ where $\pi\colon \widehat{G/G_\xi}\to D_3$ is the equivariant projection. Since $\pi$ is injective on $G/G_\xi$, $(g_n)_*\nu\to \delta_{G_\xi}$.
\end{proof}

\begin{remark} The proof of Theorem~\ref{umf} shows that the map $\widehat{G/G_\xi}\to D_3$ is 6 to 1 on preimages of branch points, 4 to 1 on preimages of regular points because one has the choice of a connected component and a choice of an orientation. The map $\widehat{G/G_\xi}\to D_3$ is 2 to 1 on points that are not preimages of end points. We have the following commutative diagram.

\begin{center}
 \begin{tikzcd}[row sep=tiny]
\widehat{G/G_\xi} \arrow[dd] \arrow[dr]&  \\
 &\bS^1  \arrow[dl]            \\
 D_3
\end{tikzcd}

\end{center}

\end{remark}
\bibliographystyle{../../Latex/Biblio/halpha}
\bibliography{../../Latex/Biblio/biblio.bib}

\begin{thebibliography}{MNTU00}

\bibitem[BdlHV08]{MR2415834}
Bachir Bekka, Pierre de~la Harpe, and Alain Valette.
\newblock {\em Kazhdan's property ({T})}, volume~11 of {\em New Mathematical
  Monographs}.
\newblock Cambridge University Press, Cambridge, 2008.

\bibitem[BM00]{burger2000groups}
Marc Burger and Shahar Mozes.
\newblock Groups acting on trees: from local to global structure.
\newblock {\em Publications Math{\'e}matiques de l'Institut des Hautes
  {\'E}tudes Scientifiques}, 92(1):113--150, 2000.

\bibitem[BT18]{bonk2018continuum}
Mario Bonk and Huy Tran.
\newblock The continuum self-similar tree.
\newblock {\em arXiv}, 2018, 1803.09694.

\bibitem[BZ20]{basso2020topological}
Gianluca Basso and Andy Zucker.
\newblock Topological dynamics beyond polish groups, 2020, 2008.08471.

\bibitem[Cal07]{MR2327361}
Danny Calegari.
\newblock {\em Foliations and the geometry of 3-manifolds}.
\newblock Oxford Mathematical Monographs. Oxford University Press, Oxford,
  2007.

\bibitem[CF06]{MR2207794}
Danny Calegari and Michael~H. Freedman.
\newblock Distortion in transformation groups.
\newblock {\em Geom. Topol.}, 10:267--293, 2006.
\newblock With an appendix by Yves de Cornulier.

\bibitem[CG93]{MR1230383}
Lennart Carleson and Theodore~W. Gamelin.
\newblock {\em Complex dynamics}.
\newblock Universitext: Tracts in Mathematics. Springer-Verlag, New York, 1993.

\bibitem[DH20]{bertr2020non}
Bertrand Deroin and Sebastian Hurtado.
\newblock Non left-orderability of lattices in higher rank semi-simple lie
  groups.
\newblock {\em arXiv preprint}, 2020, 2008.10687.

\bibitem[DM18]{DuchesneMonod}
Bruno Duchesne and Nicolas Monod.
\newblock Group actions on dendrites and curves.
\newblock {\em Annales de l'institut Fourier}, 68(5):558--565, 2018.

\bibitem[DM19]{MR3894039}
Bruno Duchesne and Nicolas Monod.
\newblock Structural properties of dendrite groups.
\newblock {\em Trans. Amer. Math. Soc.}, 371(3):1925--1949, 2019.

\bibitem[DMW19]{MR4017014}
Bruno Duchesne, Nicolas Monod, and Phillip Wesolek.
\newblock Kaleidoscopic groups: permutation groups constructed from dendrite
  homeomorphisms.
\newblock {\em Fund. Math.}, 247(3):229--274, 2019.

\bibitem[Duc20]{duchesne2019topological}
Bruno Duchesne.
\newblock Topological properties of {W}a\.zewski dendrite groups.
\newblock {\em J. Éc. polytech. Math.}, 7:431--477, 2020.

\bibitem[ET16]{MR3417738}
David~M. Evans and Todor Tsankov.
\newblock Free actions of free groups on countable structures and property
  ({T}).
\newblock {\em Fund. Math.}, 232(1):49--63, 2016.

\bibitem[GG20]{GARCIAGARRIDO2020105417}
Víctor~J. García-Garrido.
\newblock Unveiling the fractal structure of julia sets with lagrangian
  descriptors.
\newblock {\em Communications in Nonlinear Science and Numerical Simulation},
  91:105417, 2020.

\bibitem[Ghy99]{MR1703323}
\'{E}tienne Ghys.
\newblock Actions de r\'{e}seaux sur le cercle.
\newblock {\em Invent. Math.}, 137(1):199--231, 1999.

\bibitem[Ghy01]{MR1876932}
\'{E}tienne Ghys.
\newblock Groups acting on the circle.
\newblock {\em Enseign. Math. (2)}, 47(3-4):329--407, 2001.

\bibitem[GM83]{MR708367}
M.~Gromov and V.~D.~and Milman.
\newblock A topological application of the isoperimetric inequality.
\newblock {\em Amer. J. Math.}, 105(4):843--854, 1983.

\bibitem[GM93]{MR1209913}
Lisa~R. Goldberg and John Milnor.
\newblock Fixed points of polynomial maps. {II}. {F}ixed point portraits.
\newblock {\em Ann. Sci. \'{E}cole Norm. Sup. (4)}, 26(1):51--98, 1993.

\bibitem[GM06]{MR2255499}
James Giblin and Vladimir Markovic.
\newblock Classification of continuously transitive circle groups.
\newblock {\em Geom. Topol.}, 10:1319--1346, 2006.

\bibitem[GM18]{MR3821724}
Eli Glasner and Michael Megrelishvili.
\newblock More on tame dynamical systems.
\newblock In {\em Ergodic theory and dynamical systems in their interactions
  with arithmetics and combinatorics}, volume 2213 of {\em Lecture Notes in
  Math.}, pages 351--392. Springer, Cham, 2018.

\bibitem[KPT05]{MR2140630}
A.~S. Kechris, V.~G. Pestov, and S.~Todorcevic.
\newblock Fra\"{i}ss\'{e} limits, {R}amsey theory, and topological dynamics of
  automorphism groups.
\newblock {\em Geom. Funct. Anal.}, 15(1):106--189, 2005.

\bibitem[Kwi18]{MR3893292}
Aleksandra Kwiatkowska.
\newblock Universal minimal flows of generalized {W}a\.{z}ewski dendrites.
\newblock {\em J. Symb. Log.}, 83(4):1618--1632, 2018.

\bibitem[Mar00]{MR1797749}
Gregory Margulis.
\newblock Free subgroups of the homeomorphism group of the circle.
\newblock {\em C. R. Acad. Sci. Paris S\'{e}r. I Math.}, 331(9):669--674, 2000.

\bibitem[Mil06]{MR2193309}
John Milnor.
\newblock {\em Dynamics in one complex variable}, volume 160 of {\em Annals of
  Mathematics Studies}.
\newblock Princeton University Press, Princeton, NJ, third edition, 2006.

\bibitem[MNTU00]{MR1747010}
S.~Morosawa, Y.~Nishimura, M.~Taniguchi, and T.~Ueda.
\newblock {\em Holomorphic dynamics}, volume~66 of {\em Cambridge Studies in
  Advanced Mathematics}.
\newblock Cambridge University Press, Cambridge, 2000.
\newblock Translated from the 1995 Japanese original and revised by the
  authors.

\bibitem[Nad92]{Nadler}
Sam~Bernard Nadler, Jr.
\newblock {\em Continuum theory}, volume 158 of {\em Monographs and Textbooks
  in Pure and Applied Mathematics}.
\newblock Marcel Dekker, Inc., New York, 1992.
\newblock An introduction.

\bibitem[Nav02]{MR1951442}
Andr\'{e}s Navas.
\newblock Actions de groupes de {K}azhdan sur le cercle.
\newblock {\em Ann. Sci. \'{E}cole Norm. Sup. (4)}, 35(5):749--758, 2002.

\bibitem[Nav10]{MR2602845}
Andr\'{e}s Navas.
\newblock A finitely generated, locally indicable group with no faithful action
  by {$C^1$} diffeomorphisms of the interval.
\newblock {\em Geom. Topol.}, 14(1):573--584, 2010.

\bibitem[Nav11]{MR2809110}
Andr\'{e}s Navas.
\newblock {\em Groups of circle diffeomorphisms}.
\newblock Chicago Lectures in Mathematics. University of Chicago Press,
  Chicago, IL, spanish edition, 2011.

\bibitem[Nav18]{MR3966841}
Andr\'{e}s Navas.
\newblock Group actions on 1-manifolds: a list of very concrete open questions.
\newblock In {\em Proceedings of the {I}nternational {C}ongress of
  {M}athematicians---{R}io de {J}aneiro 2018. {V}ol. {III}. {I}nvited
  lectures}, pages 2035--2062. World Sci. Publ., Hackensack, NJ, 2018.

\bibitem[Pes98]{MR1608494}
Vladimir~G. Pestov.
\newblock On free actions, minimal flows, and a problem by {E}llis.
\newblock {\em Trans. Amer. Math. Soc.}, 350(10):4149--4165, 1998.

\bibitem[Ros09]{MR2503307}
Christian Rosendal.
\newblock A topological version of the {B}ergman property.
\newblock {\em Forum Math.}, 21(2):299--332, 2009.

\bibitem[Sok15]{MR3436366}
Miodrag Soki\'{c}.
\newblock Semilattices and the {R}amsey property.
\newblock {\em J. Symb. Log.}, 80(4):1236--1259, 2015.

\bibitem[Thu09]{MR2508255}
William~P. Thurston.
\newblock On the geometry and dynamics of iterated rational maps.
\newblock In {\em Complex dynamics}, pages 3--137. A K Peters, Wellesley, MA,
  2009.
\newblock Edited by Dierk Schleicher and Nikita Selinger and with an appendix
  by Schleicher.

\bibitem[Tsa12]{MR2929072}
Todor Tsankov.
\newblock Unitary representations of oligomorphic groups.
\newblock {\em Geom. Funct. Anal.}, 22(2):528--555, 2012.

\bibitem[Wa{\.z}23]{Wazewski23}
Tadeusz Wa{\.z}ewski.
\newblock {Sur les courbes de Jordan ne renfermant aucune courbe simple
  ferm\'ee de Jordan}.
\newblock {\em {Ann. Soc. Polon. Math.}}, 2:49--170, 1923.

\bibitem[Wit94]{MR1198459}
Dave Witte.
\newblock Arithmetic groups of higher {${\mathbf{ Q}}$}-rank cannot act on
  {$1$}-manifolds.
\newblock {\em Proc. Amer. Math. Soc.}, 122(2):333--340, 1994.

\bibitem[Zuc20]{zucker_2020}
Andy Zucker.
\newblock Maximally highly proximal flows.
\newblock {\em Ergodic Theory and Dynamical Systems}, page 1–21, 2020.

\end{thebibliography}
\end{document}